\documentstyle[french]{article}
\input amssym.def
\input amssym.tex

\newtheorem{theorem}{{Theorem}}
\newtheorem{theoreme}{{Th\' eor\`eme}} 
\newtheorem{corollaire}[theoreme]{{Corollaire}} 

\newtheorem{theoreme.algebrique}[theorem]{{Th\' eor\`eme alg\'ebrique}}%[section]
\newtheorem{theoreme.geometrique}[theorem]{{Th\' eor\`eme g\'eom\'etrique}}%[section]

\newtheorem{proposition}{{Proposition}} 
\newtheorem{lemme}[proposition]{{Lemme}}

\newtheorem{proposition.fondamentale} [theorem]{{Proposition fondamentale}}%[section]
\newtheorem{remarque}[theorem]{{Remarque}}%[section]
\newtheorem{remarque.terminologique}[theorem]{{Remarque terminologique}}%[section]
\newtheorem{quest.algeb}[theorem]{{Question alg\'ebrique}}%[section]
\newtheorem{conj.geom}[theorem]{{Conjecture g\'eom\'etrique}}%[section]
\newenvironment{preuve}{{\it Preuve.}}{\hfill$\diamondsuit$\medskip}
\newenvironment{esqpreuve}{{\it Esquise de la preuve.}}{\hfill$\diamondsuit$\medskip}

\begin{document}

%¥¥¥macros

%¥¥¥titre

\title{Quelques remarques sur les actions analytiques des r\'eseaux 
des groupes de Lie de rang sup\'erieur {\footnote {Version l\'eg\`erement 
d\'etaill\'ee d'une note soumise aux C.R.A.S.}}}
\author{Abdelghani Zeghib}
\date{\today}
\maketitle

\noindent
{\small{{\bf R\'esum\'e.  }
%\begin{abstract}  
Toutes les actions consid\'er\'ees ici sont {\it analytiques} 
(r\'eelles).
Soit $\Gamma$ 
un sous-groupe
d'indice fini de $SL(n, {\Bbb Z})$.
Nous 
montrons, 
 en particulier, 
 la rigidit\'e homotopique (globale) de son  action affine usuelle sur le
tore $T^n$ ($n\geq 3$), 
 ainsi que celle de son action projective usuelle sur la sph\`ere 
 $S^{n-1}$ (pour
$n \geq 4$). 
%Nous montrons aussi la rigidit\'e locale du produit de l'action
 %usuelle
%de $\Gamma$ sur $T^n$ par son action triviale sur une vari\'et\'e 
%compacte $N$ quelconque. Les d\'etails ainsi que des r\'esultats compl\'ementaires 
%para\^{\i}tront ult\'erieurement.
} }   \\ 
 
%\end{abstract}

%\noindent
%{\large {\bf Some remarks on analytic actions of lattices of Lie groups of higher rank}}  \\

\noindent
{\small{{\bf Abstract.} All the actions considered here are (real) 
{\it analytic}. Let 
$\Gamma$ be a subgroup of finite index of $SL(n, {\Bbb Z})$. We prove, in particular, 
the (global) homotopical  rigidity, for both  its standard affine action on the torus 
$T^n$ ($n \geq 3$), 
%s well as for 
and its standard projective action 
on the sphere $S^{n-1}$ ($n \geq 4$). \\
% the  homotopical 

 }}

\medskip
\noindent
{\bf Abridged English version} \\

 Let $\rho_0$ denote the standard affine action of 
$\Gamma$, a subgroup of finite index of $SL(n, {\Bbb Z})$, on the torus $T^n$.
In order to prove that it is locally rigid, it is a natural idea to 
show that every $C^1$ nearby action $\rho$ is ``linearizable''. From a result on 
persistance of fixed points of \cite{sto}, and a local linearizability
result of  \cite{C-G} (in the  analytic case), we get a local linearization 
of $\rho$ around some fixed point. The main ingredient is then, to introduce 
a {\it Siegel   neighborhood}  of the fixed point, which is a maximal
{\it invariant} open set on which the action is linear. We show that the action
on the Siegel neighborhood is $C^\omega$-conjugate to the  standard action on a 
punctured torus, i.e. $\rho_0$ restricted to the torus with some rational points removed.
%It may  happen that holes actually 
Holes may exist,  in general, as 
%this is the case 
in the case of   Katok-Lewis
 examples
which are constructed by a blowing up process.
But in the  torus case, a homological consideration leads to the existence
 of (at least) a periodic point lying in a hole. Its Siegel  neighborhood 
%of such a
%periodic point 
must cut that of our initial fixed point (because we are on a torus), contradicting 
the maximality of these neighborhoods. Therefore the Siegel  neighborhood is 
the whole torus, that is
$\rho$ is conjugate to $\rho_0$. This idea may be adapted to handle
% other rigidity situations as follows.
the following situations.

\begin{theorem} Let 
$\Gamma$ be a subgroup of finite index of $SL(n, {\Bbb Z})$, and $\rho_0$ its 
   standard action on 
$T^n$ ($n \geq 3$). 

i) The product of $\rho_0$ by the trivial action of $\Gamma$
on any compact (real) analytic manifold $N$ is $C^1$ locally rigid 
among $C^\omega$-actions, i.e. it is $C^\omega$-conjugate to every 
analytic $\Gamma$-action on $T^n \times N$, which is $C^1$-close to it.
%that is any analytic which  $C^1$-close 
%to it, is $C^\omega$-conjugate to it.

ii) The orbit $Diff^\omega(T^n).\rho_0$, i.e. the space of conjugates of
$\rho_0$ in the space $Rep(\Gamma \to Diff^\omega(T^n))$
of analytic actions of $\Gamma$ on $T^n$, is closed-open for  the 
$C^1$ topology, and closed for  the $C^0$ topology.

iii) Any faithful analytic action of $\Gamma$ on $T^n$ preserving a non-atomic measure, and 
having a fixed point in the support of this measure, is
% finite, or 
(up to 
an automorphism) $C^\omega$-conjugate to 
$\rho_0$.

\end{theorem}

\begin{theorem} Let 
$\Gamma$ be a subgroup of finite index of $SL(n+1, {\Bbb Z})$, $n \geq 3$, and 
$\rho_0$ its standard projective action on $S^n$.
The orbit $Diff^\omega(S^n).\rho_0$ 
%i.e. the space of conjugates of
%$\rho_0$ in the space $Rep(\Gamma \to Diff^\omega(S^n))$
%of analytic actions of $\Gamma$ on $T^n$, 
in $Rep(\Gamma \to Diff^\omega(S^n))$
is closed-open for  the 
$C^1$ topology and closed for  the $C^0$ topology. \\

\end{theorem}

%\begin{theorem }Let 
%$\Gamma$ be a subgroup of finite index of $SL(n, {\Bbb Z})$

%\end{theorem}

%\begin{abstract}  
 
%\end{abstract}

-----------------------------------------------------------------------

%\section{Introduction} 

\medskip
Pour montrer que l'action usuelle d'un sous-groupe d'indice fini 
de $SL(n,  {\Bbb Z})$ ($n \geq 3$) sur le tore $T^n$ est localement rigide, 
une id\'ee naturelle consiste \`a montrer que toute action proche est 
``lin\'earisable''. Dans cette note,  nous montrons comment faire 
marcher cette approche, et aussi la g\'en\'eraliser \`a  d'autres
 situations. Les preuves elles m\^emes, hormis 
un r\'esultat de lin\'earisablit\'e 
locale de \cite{C-G}, sont lin\'eaires, i.e. d\'ecoulent de propri\'et\'es 
de l'action usuelle.
% elle m\^eme.
 Les d\'etails ainsi que des r\'esultats compl\'ementaires 
para\^{\i}tront ult\'erieurement.

%Un exemple de r\'esultat qu'on
%retrouve ici est la rigidit\'e local 
% de l'action usuelle de $SL(n, {\Bbb Z})$ sur le tore $T^n$.
% La preuve proc\`ede par lin\'earisation, et du coup
%la preuve elle m\^eme devient lin\'eaire
% au sens qu'elle r\'esulte de propri\'et\'es 
% de l'action lin\'eaire.

Toutes les actions consid\'er\'ees ici sont {\bf analytiques} (r\'eelles).

\paragraph{1. Lin\'earisation: du local au global.}  Soit $M$ une vari\'et\'e
 analytique
r\'eelle et $\Gamma$ un groupe agissant  analytiquement sur $M$, 
via un homomorphisme $\rho: \Gamma \to $ Diff$^\omega(M)$.
  Supposons que cette action admette un point fixe $x_0$.Nous avons  alors 
une repr\'esentation infinit\'esimale: $r: \Gamma \to GL(T_{x_0}M)$.

Supposons 
que
% L'action 
$\rho$ soit 
  analytiquement lin\'earisable au voisinage de $x_0$, c'est-\`a-dire
qu'elle
 est conjugu\'ee \`a
$r$, dans un voisinage de $x_0$. Plus pr\'ecis\'ement, cela signifie qu'il
 existe $U$,  voisinage de $0$ dans $T_{x_0}M$,  et $V$ voisinage de $
x_0$ dans $M$,  et un diff\'eomorphisme analytique $\phi: U \to V$, tels 
que, pour tout $\gamma \in \Gamma$, nous avons  l'\'egalit\'e
$\phi r(\gamma) = \rho(\gamma) \phi$, dans un certain voisinage de $x_0$
(d\'ependant de 
$\gamma$).
% Supposons dans la suite que $\rho$ est lin\'earisable au voisinage de $x_0$. 
Pour pouvoir
en profiter dynamiquement, nous aurons  besoin d'ouverts $U$ et $V$ invariants
respectivement par $r$ et $\rho$.
%, dont l'existence d\'ecoule de 
%l'affirmation suivante.

\begin{proposition} 
${}$

1) Parmi les ouverts {\em \'etoil\'es} (en 0) de $T_{x_0}M$ dans 
lesquels $\phi$ 
se prolonge analytiquement en un diff\'eomorphisme local sur son image, il existe un 
seul 
%\'el\'ement
ouvert  maximal, not\'e ${\cal E}$. Il est en particulier invariant par $r$, 
et le prolongement $\bar{\phi}$ de $\phi$ \`a
${\cal E}$ v\'erifie la semi-conjugaison globale: $\bar{\phi} r(\gamma) (u) =
 \rho(\gamma) \bar{\phi} (u)$, $ \forall \gamma \in \Gamma$
et $\forall u \in {\cal E}$.

2) Il existe un ouvert maximal parmi les ouverts connexes, contenant 0,
 $r$-invariants, et  sur lesquels
$\phi$ se prolonge en un diff\'eomorphisme
analytique local semi-conjuguant $r$ et $\rho$. (Dans la suite, nous choisirons  un tel ouvert qu'on notera
${\cal M}$, et noterons  $\Phi: {\cal M} \to M$, le prolongement de 
$\phi$).

\end{proposition}

\begin{preuve} En g\'en\'eral, il n' y a pas une notion consistante de
domaine maximal d'extension d'une application analytique; sauf en dimension
(r\'eelle)  1, auquel cas, on peut parler d'intervalle maximal de prolongement 
analytique. Dans notre cas ici, nous consid\`erons  le
prolongement maximal de $\phi$ le long des demi-rayons issus de $0
\in T_{x_0}M$. Cela d\'etermine un ensemble \'etoil\'e 
(en 0) sur lequel $\phi$ se prolonge radialement analytiquement. Nous 
obtenons  un ouvert \'etoil\'e, en consid\'erant l'int\'erieur 
de cet ensemble.Nous nous se restreignons  ensuite aux points au voisinage desquels
le prolongement de $\phi$ est un diff\'eomorphisme analytique local.
C'est notre domaine ${\cal E}$. Tout est naturel, donc, 
${\cal E}$ est $r$-invariant, et la semi-conjugaison 
$\bar{\phi} r(\gamma) = \rho(\gamma) \bar{\phi}$ est satisfaite sur 
${\cal E}$. (Plus analytiquement, soit  $u \in T_{x_0}M$, si l'application 
$t \to \phi(tu)$ se prolonge \`a 
$[0, T[$, alors, l'application $t \to
\phi(t (r(\gamma)))$ se prolonge \'egalement 
dans $[0, T[$
par la formule $\phi(t (r(\gamma))) = \rho (\gamma) (\phi(tu))$...).

Enfin, comme il y a un ouvert connexe $r$-invariant satisfaisant la 
semi-conjugaison ci-dessus, il y en a un (ouvert connexe) maximal
(parmi tous les ouverts connexes, non n\'ecessairement \'etoil\'es).

\end{preuve}

%\begin{remarque}

%\end{remarque}

\paragraph{2. Exemples.} En d\'epit de son ``\'evidence'',  la 
 proposition ci-dessus ne semble pas 
exister dans la litt\'erature, o\`u l'on s'int\'eresse 
particuli\`erement au cas o\`u la re\-pr\'esent\-ation 
lin\'eaire $r$ 
est constitu\'ee soit de dilatations (et contractions), soit 
de transformations orthogonales. 
%d'homoth\'etie, soit 
et de plus, en g\'en\'eral, $\Gamma = {\Bbb Z}$.
(On parle alors de      disque, ou domaine, de Siegel...).
Il est vrai que  c'est dans ces cas que $\Phi: {\cal M} \to M$
jouit d'int\'eressantes propri\'et\'es topologiques.
Dans le cas g\'en\'eral, $\Phi$
n'est pas n\'ecessairement un rev\^etement sur son image.
On peut par exemple s'amuser \`a consid\'erer l'exemple 
de $\Gamma= {\Bbb Z}$ agissant sur le tore $
{T}^2 = S^1 \times S^1$, par un diff\'eomorphisme 
$f = g \small{\circ} A$, o\`u $A$ est 
un automorphisme hyperbolique de 
${T}^2$ et $g = (g_1, g_2)$, avec $g_1$ et $g_2$
diff\'eomorphismes analytiques de $S^1$ fixant  
$0$. Pour un choix g\'en\'erique des d\'eriv\'ees
$g_1^\prime(0)$ et $g_2^\prime(0)$, la matrice d\'eriv\'ee 
$D_{(0,0)}f$ n'aura pas de r\'esonances
(ce qui entra\^{\i}ne  en particulier 
que $\det(D_{(0,0)}f) \neq \pm 1$), et ainsi 
d'apr\`es 
le th\'eor\`eme de lin\'earisation de Sternberg, $f$ est lin\'earisable 
au voisinage de $(0,0)$. Il est 
%certainement 
exceptionnel que 
$\Phi: {\cal M} \to M$ soit un 
rev\^etement sur son image. On peut montrer 
qu'il est  impossible que 
la connexion plate de ${\cal M}$ ($\subset {\Bbb R}^2$)
descende par $\Phi$ 
% \`a $T^2$ 
(ou en d'autres termes
que les identifications de points de ${\cal M}$ ayant 
une m\^eme
image par 
$\Phi$, se fassent  \`a l'aide d'applications affines), et ce \`a
cause du fait que $\det(D_{(0,0)}f) \neq \pm 1$.

D'autres exemples s'obtiennent en consid\'erant l'action 
\`a gauche d'un 
r\'eseau $\Gamma$ d'un groupe de Lie $G$ sur le quotient 
$G / \Gamma$. Cette action fixe le point correspondant \`a
l'\'el\'ement 
neutre de $G$, et y  est (localement)  lin\'earisable. 
Ici, $\Phi$ est la restriction 
de l'application exponentielle 
${\cal G} \to G/ \Gamma$, \`a ${\cal M}$, qui  est un ouvert connexe 
%de l'alg\`ebre de Lie ${\cal G}$, 
  $Ad(\Gamma)$-invariant, et sur lequel l'application 
exponentielle est un diff\'eomorphisme local. M\^eme dans le cas 
de $G = SL(n, {\Bbb R})$, $\Phi$
 pr\'esente assez de pathologies topologiques.
 
\paragraph{3. Un cas rigide.} Les r\'eseaux irr\'eductibles des groupes 
de Lie semisimples, 
de centre fini,  sans facteur compact et de rang r\'eel 
$\geq 2$, e.g. $SL(n, {\Bbb R})$, $n\geq 3$,  sont bien connus
par leurs propri\'et\'es de rigidit\'e (et super-rigidit\'e).
Nous  n'en citerons qu'une que nous allons tout de suite 
 exploiter: un r\'esultat
de C. Cairns et E. Ghys \cite{C-G}, affirmant, avec les notations 
%et hypoth\`eses
ci-dessus (i.e. que $\Gamma$
agit analytiquement en fixant $x_0$),  que si  $\Gamma$ est un tel r\'eseau, alors 
son action est lin\'earisable au voisinage de $x_0$. 
Il est donc  int\'eressant de comprendre $\Phi$ dans ce cas. 
L'exemple pr\'ec\`edent de $\Gamma$ agissant sur 
$G/ \Gamma$, montre les limites d'une rigidit\'e esp\'er\'ee
pour $\Phi$. La raison, dans ce cas, r\'eside, probablement dans 
la ``relative pauvret\'e dynamique'' de la repr\'esentation lin\'eaire
$r$, qui n'est rien d'autre que la repr\'esentation 
$Ad$ de $\Gamma$; par exemple, $Ad(\Gamma)$ agit proprement sur 
un ouvert (non-vide) de ${\cal G}$.

Dans la suite, nous traiterons le cas o\`u $\Gamma$ est un sous-groupe
de $SL(n, {\Bbb R})$ agissant sur une vari\'et\'e $M$ de dimension ($n \geq 3$).

Nous nous  ram\`enons  au cas o\`u 
$r$ est fid\`ele (nous  utilisons pour cela  un th\'eor\`eme 
de Margulis affirmant que le noyau de $r$ est soit fini, 
soit d'indice fini).  
D'apr\`es la super-rigidit\'e de Margulis,  $r$ est soit  
 la repr\'esentation canonique soit  sa duale; 
nous    supposerons  pour simplifier les notations que
%$r$ est la repr\'esentation 
c'est la canonique. 
%(les deux autres possibilit\'es
%sont que $r$ soit finie ou la duale de la canonique).

%La repr\'esentation $r$ est alors soit 
%finie, soit la repr\'esentation canonique, soit sa duale. 
% Lorsque $r$ est finie, $\rho$ elle aussi est finie, ce qu'on consid\`ere ici
%comme une situation triviale.

%Nous supposerons pour simplifier les notations que
%$r$ est la repr\'esentation canonique. 

On peut commencer par 
essayer de comprendre ${\cal M}$.
Son compl\'ementaire ${\cal C}= {\Bbb R}^n - {\cal M}$ est un ferm\'e
invariant par l'action de $\Gamma$ sur ${\Bbb R}^n$. De tels
 ensembles ont \'et\'e
intensivement \'etudi\'es dans la litt\'erature; mais tous les
 r\'esultats les concernant 
peuvent se d\'eduire, apr\`es manipulation alg\'ebrique, du
 Th\'eor\`eme
de M. Ratner, r\'esolvant la conjecture de Raghunathan \cite{Rat}. Pour voir que ce th\'eor\`eme 
s'applique bien, 
on identifie  ${\Bbb R}^n-\{0\}$ \`a l'espace homog\`ene 
$SL(n, {\Bbb R}) /H$, o\`u
$H$ est le stabilisateur d'un certain point. Une partie 
ferm\'ee  ${\cal C}$ de ${\Bbb R}^n-\{0\}$, $\Gamma$-invariante, 
s'identifie aussi  \`a
une partie ferm\'ee ${\cal C}^\prime$
 de $SL(n, {\Bbb R})$, $\Gamma$-invariante \`a
gauche et 
$H$-invariante \`a
droite, i.e. ${\cal C}^\prime= \Gamma. {\cal C}^\prime. H$; et par cons\'equent, 
elle s'identifie \`a une partie ferm\'ee de
$\Gamma \setminus SL(n, {\Bbb R})$, invariante par 
$H$ (agissant \`a droite). Le th\'eor\`eme s'applique car
$\Gamma$ est un r\'eseau, et $H$ est engendr\'e
par ses \'el\'ements unipotents.

Le th\'eor\`eme de Ratner affirme que pour $u \in {\Bbb R}^n -\{0\}$, 
 la composante connexe de $u$ dans l'adh\'erence 
$\overline{\Gamma. u}$ de son orbite, est de la forme $G.u$, o\`u
$G$ est un sous-groupe de Lie connexe,
contenant le stabilisateur de $u$ (qui est un 
conjugu\'e de $H$)  et 
tel que $\Gamma \cap G$ soit un r\'eseau de 
$G$. 

On peut en d\'eduire (mais c'est aussi faisable \`a l'aide 
d'outils plus \'el\'ementaires, dans ce cas pr\'ecis) que 
si $\Gamma$ 
n'est pas, 
%\`a conjugaison pr\`es, comensurable \`a
\`a automorphisme pr\`es, un sous-groupe d'indice fini de
$SL(n, {\Bbb Z})$, alors, 
$\Gamma$ agit minimalement sur 
${\Bbb R}^n-\{0\}$, et en particulier ${\cal C}=
\emptyset$. Si $\Gamma$ est 
\`a automorphisme pr\`es un sous-groupe d'indice fini de
$SL(n, {\Bbb Z})$, alors  l'orbite 
d'un point est, soit dense, soit 
discr\`ete, auquel cas, ce point   est de la
forme
$\lambda u$, o\`u $u$ est rationnel (et $\lambda \in {\Bbb R}$).
En g\'en\'eral, nous avons des parties ferm\'ees invariantes propres
de la forme
$\Gamma. (F.u)$ ($= F.(\Gamma.u)$), o\`u $u$ est rationnel et
%  $0 < \alpha \leq \beta \leq \infty$ 
$F$ est un ferm\'e de $[0, \infty [$ ne contenant pas 
0 (la notation $F.u$ d\'esigne 
l'ensemble  $\{ tu, t \in F \}$). Tout ferm\'e invariant propre
est r\'eunion finie de tels ferm\'es \'el\'ementaires.

%Une application 

%(ayant un nombre fini de composante 

\begin{theoreme} Soit $\Gamma$ un r\'eseau de 
$SL(n, {\Bbb R}), n>2$,  agissant analytiquement en fixant un point $x_0$ sur 
une vari\'et\'e
$M$ de dimension $n$. Alors, ${\cal M} $ est unique, on 
appellera  $\Phi({\cal M})$ l'{\bf ouvert de Siegel} en $x_0$.

Le prolongement maximal 
%de la lin\'earisante
$\Phi: {\cal M} \to M$ v\'erifie:

i) $\Phi$ est un  rev\^etement sur son image.

ii) La connexion affine plate sur ${\cal M}$ descend par $\Phi$
en une connexion plate $\nabla$ sur $\Phi({\cal M})$.

iii)
 $\nabla$ ne peut pas se prolonger (m\^eme localement) en dehors de 
l'ouvert de Siegel.

 %$\nabla$ ne se prolonge pas (m\^eme localement) en dehors de 
%$\Phi({\cal M})$. 

iv) $\nabla$ est (localement)  unique, au sens que si $\Gamma$ pr\'eserve 
une connexion d\'efinie sur un ouvert invariant contenu dans
 l'ouvert de Siegel, alors, c'est la restriction de $\nabla$.  \\
%De plus, on a l'une des trois possibilit\'es suivantes:
De plus, il y a deux situations possibles: \\
\underline{Cas dissipatif: } On y distingue deux cas.

1) $\Gamma$ n'est pas 
\`a automorphisme pr\`es un sous-groupe d'indice fini de
$SL(n, {\Bbb Z})$, alors ${\cal M} = {\Bbb R}^n$,  $\Phi$
est un diff\'eo\-morphisme (global) sur son
image, et $\nabla$ est compl\`ete.

% La boule topologique 
%$\Phi({\cal M})$ sera appel\'ee dans ce cas 
%boule de Siegel en $x_0$.

2)  $\Gamma$ est \`a automorphisme pr\`es un sous-groupe d'indice fini de
 $SL(n, {\Bbb Z})$, 
${\cal M}\subset {\Bbb R}^n$ est un ouvert $\Gamma$-invariant
comme d\'ecrit ci-dessus, et 
  $\Phi$ est un diff\'eo\-morphisme sur son image 
$\Phi({\cal M})$. \\
\underline{Cas conservatif: }  $\Gamma$ est
\`a automorphisme pr\`es un sous-groupe d'indice fini de
$SL(n, {\Bbb Z})$, et 
$\Phi$ n'est pas injective. Dans ce cas, 
$\Phi$  transite  \`a travers un diff\'eomorphisme
$ \Phi^*: (T^n)^* \to \Phi({\cal M})$, o\`u
$(T^n)^*$ est le tore $T^n$ priv\'e d'un nombre fini 
de  points rationnels, et 
$\Phi^*$ respecte les actions et les connexions.
% L'image de $\Phi$ sera appel\'ee tore (trou\'e) de Siegel.  

La  possibilit\'e conservative se produit  exactement 
lorsque
l'une ou l'autre des deux conditions suivantes est satisfaite:

1)  $\Gamma$  pr\'eserve une mesure 
finie sans atomes dont le support contient $x_0$. Cette mesure est alors
la mesure standard sur le tore trou\'e.

2) Il existe un \'el\'ement $\gamma \in \Gamma$  
d'ordre infini tel que  l'ensemble 
des points non-errants de $\rho(\gamma)$
%$\rho(\gamma)$  pr\'eserve une mesure 
% dont le support
 contient un voisinage de $x_0$.

%finie sans atomes dont le support contient $x_0$. Cette mesure est alors
%la mesure standard sur le tore trou\'e.

\end{theoreme}

\begin{esqpreuve}
D'apr\`es la description ci-dessus des ensembles ferm\'es 
invariants par  $\Gamma$ agissant sur ${\Bbb R}^n$, nous tirons  en
particulier que tout ouvert invariant est connexe. En particulier, s'il y
a deux ouverts maximaux sur lesquels $\phi$
se prolonge analytiquement, alors leur 
intersection est connexe, et par suite, $\phi$ se prolonge \`a leur r\'eunion, 
donc ${\cal M}$ est unique.

Les propri\'et\'es de $\Phi$ d\'ecoulent de la rigidit\'e
de l'ensemble: $${\cal R} = \{ (x, y) \in 
{\cal M} \times {\cal M} \subset {\Bbb R}^n \times
{\Bbb R}^n / \Phi(x) = \Phi(y) \}$$
C'est une sous-vari\'et\'e analytique de 
${\cal M} \times {\cal M}$, localement ferm\'ee dans 
${\Bbb R}^n \times {\Bbb R}^n$, invariante par 
l'action diagonale de $\Gamma$ 
sur ${\Bbb R}^n \times {\Bbb R}^n$. Pour la d\'ecrire, nous consid\`erons  
l'action diagonale de 
$SL(n, {\Bbb R})$ sur ${\Bbb R}^n \times {\Bbb R}^n-\{(0,0)\}$. Cette action 
admet une orbite ouverte 
${\cal O} = \{ (u, v) / {\Bbb R}u \neq {\Bbb R}v \}$, et 
des orbites d\'eg\'en\'er\'ees de 
 la forme $O_\alpha = \{(u, \alpha u), u \in {\Bbb R}^n \}$, $\alpha$
\'etant un r\'eel. (On suppose $n >2$).

\'Evidemment, ${\cal R}$ qui est une relation d'\'equivalence, contient 
la diagonale
${\cal O}_1$. 
L'\'egalit\'e ${\cal R} = {\cal O}_1$ 
signifie que $\Phi$ est injective.
On se convainc facilement que si 
${\cal R} \neq {\cal O}_1$,  alors 
${\cal R}$ contient des points de l'orbite ouverte ${\cal O}$.

Choisissons un point $(e_1, e_2)$ de 
${\cal O}$. Son stabilisateur $H$  (dans 
$SL(n, {\Bbb R})$) est engendr\'e par ses \'el\'ements unipotents; 
pour $n=3$, $H$ est unipotent. Pour \'etudier l'adh\'erence de 
$\Gamma. (e_1,e_2)$, nous aurons  besoin de comprendre les sous-groupes connexes 
contenant $H$.

\begin{lemme} Notons  $P$ le plan 
${\Bbb R}e_1 \bigoplus {\Bbb R}e_2$, et $S_P$
le sous-groupe de $SL(n, {\Bbb R})$ pr\'eservant 
$P$; nous avons   une projection 
$\pi: S_P \to GL(P) \simeq GL(2, {\Bbb R})$ ($H$ 
s'identifie \`a
$\pi^{-1}(1)$).

 Soit $G$ un groupe de Lie connexe contenant 
$H$ et diff\'erent de $SL(2, {\Bbb R})$,  alors l'une des 
deux possibilit\'es suivantes se pr\'esente:

i) $G \subset S_P$; plus exactement $G = 
\pi^{-1}(L)$, o\`u 
$L$ est un sous-groupe de $GL(2, {\Bbb R})$.

ii) Il existe $e \in P$, et $G = S_e$, ou 
$G= S_{ {\Bbb R}e}$, o\`u $S_{\{ \}}$
d\'esigne le stabilisateur.

%$ ou bien $H$ laisse invariant le plan 
%${\Bbb R}e_1 \bigoplus {\bf R}e_2$, ou bien 
%$H$ laisse invariante une droite contenue dans ce dernier plan.

\end{lemme}

Dans le cas de $\Gamma =
SL(n, {\Bbb Z})$, il existe $G$ avec 
$H \subset G \subset S_P$, et $G \cap \Gamma$, un r\'eseau de 
$G$, si et seulement si $P$ est un 2-plan
rationnel (i.e. engendr\'e par deux vecteurs rationnels). Pour $\Gamma$
quelconque,
% Il est en g\'en\'eral  vrai que
 l'ensemble des 2-plans $P$ avec un groupe $G$ comme pr\'ec\'edemment, 
v\'erifiant que $G \cap \Gamma$ est un r\'eseau de $G$, est d\'enombrable.

On peut donc se restreindre \`a l'\'etude du cas
o\`u l'adh\'erence de l'orbite 
$\Gamma.(e_1, e_2)$ est d\'etermin\'ee par un groupe $G$
\'egal \`a $S_e$ ou $S_{{\Bbb R}e }$. En fait $G= S_e$ 
 car
  $S_{{\Bbb R}e}$
 n'est 
pas unimodulaire.  Donc pour l'action de $\Gamma$
sur ${\Bbb R}^n$, l'orbite de $e$ est discr\`ete; il en d\'ecoule
que 
$\Gamma$ est 
\`a automorphisme pr\`es un sous-groupe d'indice fini de
$SL(n, {\Bbb Z})$.

L'\'ev\'enement,   $e$ est colin\'eaire \`a $e_1$
ou \`a  $e_2$, correspond \`a un nombre d\'enombrable
de cas; nous pouvons  donc supposer qu'il n'a pas lieu.

Maintenant $S_e$ agit transitivement sur 
${\Bbb R}^n -{\Bbb R}e$; en particulier tous  
les points $(x, y)$ d'un voisinage de  $(e_1,e_2)$
dans ${\cal R}$, 
s'\'ecrivent: $(x, y)= (A(e_1), A(e_2))$, 
avec $A \in S_e$.

Comme $e$ est d\'efini \`a
facteur pr\`es, et appartient \`a 
$P$ ($= {\Bbb R}e_1 \bigoplus {\Bbb R}e_2$), nous pouvons  \'ecrire 
$e_2= \alpha e_1 +e$, pour un certain $\alpha$. Donc $y = A(e_2) = \alpha
A(e_1) + A(e)
= \alpha  x + e$ (car $A \in S_e$). 

En d'autres termes, un voisinage 
de $(e_1, e_2)$
dans ${\cal R}$ co\"{\i}ncide 
avec le graphe de l'homoth\'etie-translation $x \to 
\alpha x +e$. Notons $f$ cette transformation.
Alors, dans
un voisinage de $e_1$, nous avons   $\Phi \small{\circ} f = \Phi$. Mais comme $f$
 est d\'efinie
partout, l'\'egalit\'e se prolonge \`a tout ${\cal M}$, 
qui doit \^etre invariant par $f$.

Il en va de m\^eme pour les  transformations  $h= AfA^{-1}$,
qui sont de la forme $x \to \alpha x + A(e)$,  o\`u 
$A$ parcourt $\Gamma$. Sachant que 
$\Gamma$
est \`a automorphisme pr\`es un sous-groupe d'indice fini de
$SL(n, {\Bbb Z})$, nous en d\'eduisons  
que les parties translationnelles $A(e)$ de ces transformations 
forment un r\'eseau de ${\Bbb R}^n$, et en particulier, si 
$\alpha \neq 1$, alors il existe un commutateur 
de ces transformations qui est une translation non triviale; et 
par suite, il existe un r\'eseau de translations $g$
v\'erifiant 
$\Phi \small{\circ} g = \Phi$.
Si $\alpha \neq \pm 1$, il existerait trop de transformations
$h$ satisfaisant $\Phi \small{\circ} h= \Phi$, contredisant 
le fait que $\Phi$ est un diff\'eomorphisme local. Ceci montre que $\Phi$
transite \`a travers un diff\'eomorphisme 
$\Phi^*$ comme \'enonc\'e.

La connexion plate descend \`a  l'ouvert de Siegel $\Phi{\cal M})$
car le groupe de Galois du rev\^etement $\Phi$ agit affinement sur 
${\cal M}$. Le fait qu'elle ne peut pas se prolonger \`a
un ouvert 
strictement plus grand que l'ouvert de Siegel, provient du fait que
cette connexion est d\'ej\`a ``pratiquement compl\`ete''.
\end{esqpreuve}

\begin{corollaire}  Si les ouverts de Siegel de deux points  fixes de 
$\Gamma$ s'intersectent, alors ils sont identiques.

\end{corollaire}

\paragraph{4. Dimension $>n$.} La discussion pr\'ec\'edente se g\'en\'eralise 
partiellement lorsque $\Gamma$ agit sur une vari\'et\'e $M$ de 
dimension $n +p$, telle  que la repr\'esentation $r$ de $\Gamma$ dans
${\Bbb R}^n \times {\Bbb R}^p$ soit  le produit de la repr\'esentation 
canonique dans ${\Bbb R}^n$ par la repr\'esentation triviale dans 
${\Bbb R}^p$. 

Par exemple, on peut se   ramener essentiellement \`a
cette situation, lorsque $p <n$. 

Les transformations 
(partielles) $f$ de ${\Bbb R}^{n+p}$ v\'erifiant 
$\Phi \small{\circ} f = \Phi$ sont de la forme:
$$f: (u,v) \in {\Bbb R}^n \times U \subset {\Bbb R}^n \times {\Bbb R}^p \to
(\alpha(v)u + e(v), g(v))$$ o\`u $\alpha: U \to {\Bbb R}$, 
$e: U\to {\Bbb R}^n$ et $g: U \to {\Bbb R}^p$, sont des applications
analytiques 
d\'efinies sur  l'ouvert  (de ${\Bbb R}^p$)  $U$ et  $g$ est un diff\'eomorphisme
sur son image. (\'Evidemment, il peut se passer que
$\alpha(v)  \equiv  1$, et $e(v) \equiv 0$).

Notons que vu la forme de $f$, la composition, 
\`a gauche ou \`a droite, de $f$ avec des \'el\'ements de 
la forme $r(\gamma)$, est d\'efinie dans tout le domaine de d\'efinition de
$f$.

Le produit $g =  r(\gamma) f^{-1} r(\gamma)^{-1} f$, est un diff\'eomorphisme 
de ${\Bbb R}^n \times U$  de la
 forme
$g: (u, v) \to (u+ e(v)- r(\gamma) (e(v)), v)$. En particulier, 
d\`es qu'il existe un $e(v) \neq 0$, alors, il y 'en aura un r\'eseau,  $\Gamma$ 
est \`a automorphisme pr\`es un sous-groupe d'indice fini de
$SL(n, {\Bbb Z})$, et le niveau 
${\Bbb R}^n \times \{v\}$ s'envoie par $\Phi$ sur un tore 
trou\'e.  

%Notons aussi que si $R(v) \subset {\Bbb R}^n$ est le r\'eseau 
%%correspondant \`a 
%engendr\'e par les vecteurs $\{e(v) \}$, alors au voisinage
%d'un vecteur $v_0$, on a 
%$R(v) = t(v)R(v_0)$, o\`u $t(v)$ est une fonction r\'eelle analytique.

%Soit $\gamma \in Gamma$, et notons $A = r(\gamma)$, et 
%$f_A =AfA^{-1}$ (notons que vu la forme de $f$, la composition de 
%(notons que vu la forme de $f$, la composition de 

Cette discussion nous fournit en particulier le fait suivant.

\begin{corollaire} Soit $\Gamma$ un r\'eseau de
$SL(n, {\Bbb R})$ agissant analytiquement sur une vari\'et\'e $M$, en 
fixant un point $x_0$, avec une repr\'esentation infinit\'esimale, 
produit de la  repr\'esentation canonique dans ${\Bbb R}^n$ par une 
repr\'esentation triviale. Supposons que l'action pr\'eserve une mesure
finie dont le support contient un voisinage de $x_0$. Alors, 
$\Gamma$ est \`a automorphisme pr\`es un sous-groupe d'indice fini de
$SL(n, {\Bbb Z})$, et sur un  ouvert invariant contenant 
$x_0$, l'action est conjugu\'ee \`a l'action sur 
$(T^n)^* \times N$, o\`u $(T^n)^*$  est un tore 
trou\'e 
sur lequel $\Gamma$ agit de la fa\c con usuelle, et 
$N$ est une vari\'et\'e de dimension $p= $ dim$M-n$, sur lequel
$\Gamma$ agit trivialement.

\end{corollaire}

On en d\'eduit en particulier:

\begin{corollaire} Soit $\Gamma$ un r\'eseau de 
$SL(n, {\Bbb R})$, $n>2 $, agissant fid\`element analytiquement   sur  
une vari\'et\'e $M$ de dimension $<2n$, en pr\'eservant
 une mesure finie pleine (i.e.son support est \'egal \`a $M$).
%$\mu$
 %lisee. 
% non atomique.
% Si l'action n'est pas finie, et
Si  $\Gamma$ n'est pas
\`a automorphisme pr\`es un sous-groupe d'indice fini de
 $SL(n, {\Bbb Z})$,  
alors, l'action n'a  aucun  point p\'eriodique.
% $\Gamma$ n'admet pas de point p\'eriodique
%dans le support de $\mu$. Plus pr\'ecis\'ement, le stabilisateur de tout point 
%du support de $\mu$ est fini.

\end{corollaire}

\begin{remarque} {\em Dans tous les \'enonc\'es pr\'ec\'edents, $M$ n'\'etait pas 
suppos\'ee compacte! }

\end{remarque}

%\part{Applications}

%\paragraph{Rigidit\'e locale de l'action
% usuelle de $SL(n, {\Bbb Z})$
%sur le tore.}

%\paragraph{5. Application: actions sur le tore.} 
\paragraph{5. Premi\`ere application: actions sur le tore.} 

Nous esquissons dans ce qui suit une preuve 
de la rigidit\'e locale de l'action usuelle d'un sous-groupe $\Gamma$
 d'indice fini dans 
$SL(n, {\Bbb Z})$ sur le tore (parmi les actions $C^\omega$).
En fait, pour  simplifier les notations, nous supposerons
 que $\Gamma = SL(n, {\Bbb Z})$. Nous supposerons aussi que $n \geq 6$, et 
ce pour pouvoir ``scinder'' $\Gamma$ en 
r\'eseaux de rang sup\'erieur. Le cas g\'en\'eral, c'est-\`a-dire $n\geq 3$, 
demande un peu plus d'analyse.

 Soit 
$\rho$ une action $C^1$ proche de l'action standard  $\rho_0$. 
D'apr\`es un th\'eor\`eme de 
Stowe \cite{sto}, $\rho$ admet  un point fixe 
$x_0$ (proche de 0). On peut se convaincre facilement que
 l'ouvert de Siegel
en $x_0$ est un tore trou\'e,
c'est-\`a-dire qu'on ne peut pas \^etre dans le cas dissipatif du 
Th\'eor\`eme 1.
%(et non pas 
%hom\'eomorphe \`a
%un ouvert de ${\Bbb R}^n$, comme dans 
%le cas compl\`etement dissipatif). 
L'action $\rho$ sera 
conjugu\'ee \`a
$\rho_0$
si l'on d\'emontre que l'ouvert de Siegel est un vrai tore, i.e. 
qu'il est 
non-trou\'e.

Pour montrer cela, il suffit
 %essentiellement 
de montrer que $\Gamma$, ou m\^eme un 
  sous-groupe d'indice fini, poss\`ede un point fixe, disons
$x_1$, r\'esidant dans un  trou (i.e.
le  compl\'ementaire de l'ouvert de Siegel de $x_0$). En effet, alors,
puisqu'on est sur un tore, 
les deux domaines de Siegel  doivent   s'intersecter, 
contredisant le Corollaire 2.
%et donc d'apr\`es  le th\'eor\`eme 1,
% ils d\'eterminent $\gamma$  

On 
%trouvera 
cherchera un point fixe dans un trou, comme intersection 
de deux sous-vari\'et\'es,  lieux de points fixes de deux sous-groupes
 engendrant $\Gamma$, dont on
% sait 
%connait 
sait que le nombre d'intersection total
%pour des raisons homologiques, mais on sait que leur
est sup\'erieur \`a leur nombre d'intersection    dans 
l'ouvert de Siegel
en $x_0$.

 Plus pr\'ecis\'ement, soit $e_1, \ldots, e_n$ la base canonique
de ${\Bbb R}^n$, et notons $E_1$ (resp. $E_2$) les plans engendr\'es 
par $e_1, e_2, e_3$ (resp. $e_4, \ldots, e_n$). Ils sont 
(ponctuellement) fix\'es 
par $ \Gamma_1= SL(n-3, {\Bbb Z})$ et $\Gamma_2= SL(3, {\Bbb Z})$ (plong\'es
 naturellement dans 
$SL(n, {\Bbb Z})$).

Notons $F_1$ (resp. $F_2$) l'ensemble (analytique) des points fixes
de $\rho(\Gamma_1)$ (resp. $\rho(\Gamma_2)$). D'apr\`es un th\'eor\`eme 
de Stowe, lorsque $\rho$ est $C^1$
proche de $\rho_0$, ces ensembles sont des sous vari\'et\'es  (proches 
des celles qui correspondent \`a $\rho_0$). (On utilise pour appliquer le 
th\'eor\`eme
de Stowe, le th\'eor\`eme d'annulation cohomologique de Margulis). 

\'Evidemment $F_1$ et $F_2$ sont des prolongements dans les trous,  de 
$\Phi(E_1)$ et $\Phi(E_2)$ respectivement. 

Consid\'erons  un \'el\'ement 
$\gamma \in SL(n, {\Bbb Z})$ tel que 
$E_2^\prime = r(\gamma)(E_2)$ 
%$r(\gamma)(E_2)$ 
intersecte transversalement 
 %l'intersection $E_2^\prime  \cap (
$E_1 + {\Bbb Z}^n$ en  un  point de ${\Bbb R}^n-{\cal M}$ 
(cela existe car 
${\Bbb R}^n-{\cal M}$ est constitu\'e de points rationnels). 

Le nombre
d'intersection de $F_2^\prime = \rho(\gamma) (F_2) $ avec $F_1$ est le m\^eme que dans le
cas de l'action standard. Il est \'egal au cardinal de $E_2^\prime  \cap (E_1 + {\Bbb Z}^n)$
mod. ${\Bbb Z}^n$. Mais le nombre d'intersection de $F_2^\prime$ avec $F_1$
\`a l'int\'erieur de l'ouvert de Siegel
vaut le nombre d'intersection de $\Phi(E_2^\prime)$ avec $\Phi(E_1)$, qui est \'egal
au cardinal de $(E_2^\prime  \cap (E_1 + {\Bbb Z}^n)  -{\cal M})$ mod. 
${\Bbb Z}^n$. Le choix de $\gamma$ assure que ces deux nombres d'intersection sont diff\'erents, et 
par suite $F_2^\prime$ et $F_1$ s'intersectent en dehors de l'ouvert de Siegel.

On remarque maintenant que $E_1$ et $E_2$  
sont  
en fait fix\'es (ponctuellement) 
%les ensembles de 
% points fixes de groupes plus grands que 
par des groupes plus grands que $\Gamma_2$ et $\Gamma_1$ respectivement; et 
qui sont
des produits semi-directs \'evidents $\Gamma_2 \ltimes N_2$ et 
$\Gamma_1 \ltimes N_1$ , o\`u $N_1$ et $N_2$ sont unipotents.  Il en r\'esulte
que 
le groupe engendr\'e par $r(\gamma) (\Gamma_2 \ltimes N_2) r(\gamma)^{-1}$   
et $\Gamma_1 \ltimes N_1$ fixe un point dans un trou. Mais, 
\`a cause de la transversalit\'e entre  $E_2^\prime$ et $E_1$,  ce
 groupe est d'indice 
fini dans $SL(n, {\Bbb Z})$.
\hfill {$\diamondsuit$}

Les premiers r\'esultats de  rigidit\'e locale  sont dus  
\`a Hurder \cite{Hur} et Katok-Lewis \cite{K-L}  
(dans le cas lisse). La rigidit\'e locale de
l'action des sous-groupes d'indice fini de $SL(n, {\Bbb Z})$
sur $T^n$ est parue dans 
%dansd\^ue \`a 
\cite{KLZ}. L'un des travaux les plus r\'ecents sur la question est
 \cite{M-Q}, o\`u
l'on d\'emontre la rigidit\'e locale des ``actions alg\'ebriques''
faiblement hyperboliques.
 Ici, \`a  l'aide des d\'eveloppements du \S 4, on peut adapter
% cette 
%Cette
%se g\'en\'eralise s'adapte 
 l'approche ci-dessus
pour traiter la rigidit\'e locale
dans  une 
%nouvelle 
situation 
%assez  g\'en\'erale,
qui n'est pas faiblement hyperbolique, 
% celle
o\`u  $\Gamma$ agit diagonalement sur un produit $T^n \times N$, o\`u
$N$ est une vari\'et\'e compacte quelconque 
 sur laquelle $\Gamma$ agit trivialement (Voir \cite{N-T} pour des r\'esultats 
proches).

\begin{theoreme} Soit 
$\Gamma$ un 
 sous-groupe d'indice fini de $SL(n, {\Bbb Z})$, $n >2$.
 Consid\'erons 
$\rho_0$ l'action produit de l'action usuelle 
de $\Gamma$
sur $T^n$, par 
l'action triviale sur une vari\'et\'e compacte 
$N$. Alors $\rho_0$ est localement rigide sous 
perturbation $C^1$, parmi
les actions analytiques; plus pr\'ecis\'ement, toute
 action 
analytique
% $\rho$, 
$C^1$ proche de $\rho_0$, est 
analytiquement conjugu\'ee \`a
$\rho_0$.

\end{theoreme}

Pour les actions sur $T^n$, nous avons  le r\'esultat de 
non-d\'eg\'en\'erescense suivant.

\begin{theoreme} Soit $\Gamma$ un sous-groupes d'indice fini de 
$SL(n, {\Bbb Z})$, et 
$(\rho_i)$ une suite d'actions analytiquesde $\Gamma$
conjugu\'ees \`a  son action usuelle sur $T^n$. Supposons que   
cette suite converge au sens de la topologie $C^0$ vers une
 action analytique $\rho$. Alors $\rho$ 
est analytiquement conjugu\'ee \`a
l'action usuelle.

\end{theoreme}

Ce type de  r\'esultat est \'etranger \`a la th\'eorie hyperbolique, 
notamment dans le cas classique des diff\'eomorphismes d'Anosov,
puisqu'il existe des applications dites DA, d\'eriv\'ees d'Anosov.

Le th\'eor\`eme  suivant 
%d\'ecoulent des 
%r\'esument 
unifie les deux  r\'esultats
pr\'ec\'edents dans le cas de $T^n$, et rend vraisemblable 
une rigidit\'e globale des actions analytiques  
%de $\Gamma$ sur $T^n$ 
dans ce cas (sans hypoth\`ese d'hyperbolicit\'e ou de 
 pr\'eservation de volume...).

\begin{corollaire}
%[rigidit\'e homotopique]
 Soit
 $\Gamma$ un sous-groupes d'indice fini de 
$SL(n, {\Bbb Z})$, et 
$Rep(\Gamma,  Diff^\omega(T^n))$
l'espace de ses  actions analytiques   sur $T^n$. 
L'orbite 
$Diff^\omega(T^n).\rho_0$, i.e. l'espace des
 actions conjugu\'ees \`a l'action standard $\rho_0$, est   
ouvert dans 
$Rep(\Gamma,  Diff^\omega(T^n)$ au sens de la topologie
$C^1$,  et ferm\'ee au sens de la topologie 
$C^0$. 

En particulier, toute action homotope \`a l'action usuelle, au sens
de la topologie $C^1$, \`a travers des actions $C^\omega$, lui est 
$C^\omega$ conjugu\'ee.

\end{corollaire}

Nous ne savons pas d\'emontrer 
que l'orbite $Diff^\omega(T^n).\rho_0$ est ouverte dans 
$Rep(\Gamma,  Diff^\omega(T^n))$, au sens de la
topologie $C^0$, pour la
simple raison que nous ne disposons pas d'un 
r\'esultat de persistance
de point fixe de $\Gamma$ sous perturbation $C^0$. (Une telle 
persistance para\^{\i}t vraisemblable dans ce contexte pr\'ecis).

Enfin,  nous avons ce r\'esultat global, qui 
r\'eduit (essentiellement) la 
rigidit\'e globale \`a l'existence de points
 p\'eriodiques.

\begin{theoreme} Toute action fid\`ele analytique d'un 
sous-groupe d'indice fini de $SL(n, {\Bbb Z})$  
sur $T^n$ ayant un point 
%p\'eriodique 
fixe et pr\'eservant une mesure finie 
non-atomique dont le support contient le point fixe est 
%(\`a automorphisme pr\`es), soit finie, soit 
analytiquement conjugu\'ee \`a l'action usuelle.

\end{theoreme}

%\paragraph{Un exemple d'application.}

\paragraph{6.  Deuxi\`eme application: actions  sur la sph\`ere.} 
%\paragraph{6. R\'esultats sur les actions projectives.}
Le groupe de
%Soit $\Gamma$ un r\'eseau 
Lie $SL(n+1, {\Bbb R})$ agit projectivement sur la sph\`ere
$S^n$, ce qui donne par restriction  une action de ses r\'eseaux.
Une super-rigidit\'e dans ce contexte consiste \`a se demander 
si r\'eciproquement, 
une action d'un r\'eseau $\Gamma$ de $SL(n+1, {\Bbb R})$ 
se prolonge en une action de $SL(n+1, {\Bbb R})$ (qui sera par suite 
n\'ecessairement l'action projective
 usuelle, \`a automorphisme
 pr\`es).  La rigidit\'e locale est une version locale de cette 
question.  Elle a \'et\'e d\'emontr\'ee, mais seulement pour 
les r\'eseaux co-compacts, d'abord par Kanai \cite{Kan} avec 
une restriction (technique) sur la dimension, et ensuite
dans le cas g\'en\'eral par Katok-Spatzier \cite{K-S}
(qui d\'emontrent en fait la rigidit\'e locale des r\'eseaux co-compacts
de rang $\geq 2$ agissant sur des bords).

Ici, nous consid\'erons des r\'eseaux (non co-compacts)
de $SL(n+1, {\Bbb R})$
qui sont  des sous-groupes d'indice fini 
dans $SL(n+1, {\Bbb Z})$ (\`a automorphisme  pr\`es). Leur avantage
 ici, est qu'il contiennent des r\'eseaux de 
$SL(n, {\Bbb R})$, qui admettent  donc des  points fixes, auxquels
nous pouvons  ainsi appliquer Th\'eor\`eme 1.
Nous avons le r\'esultat (global) suivant.

\begin{theoreme}
\label{projective.rigid}
 Soit $\Gamma$ un sous-groupe 
  d'indice fini 
de $SL(n+1, {\Bbb Z})$, agissant fid\`element 
analytiquement sur une vari\'et\'e compacte 
$M$ de dimension $n \geq 3$, tel que $\Gamma \cap SL(n, {\Bbb Z})$ 
admet un point fixe. Alors, \`a automorphisme pr\`es, 
%il s'agit, soit d'une action finie, soit  de 
 c'est l'action usuelle de $\Gamma$ sur la sph\`ere $S^n$
ou sur l'espace projectif ${\Bbb R}P^n$.

\end{theoreme}

%{\it  Quelques \'etapes de la preuve du Th\'eor\`eme
% \ref{projective.rigid}.}

\begin{esqpreuve} Nous  supposons pour simplifier les notations que $\Gamma = 
SL(n+1, {\Bbb Z})$.  
Ses \'el\'ements de la forme 
%On voit $SL(n, {\Bbb R})$ 
%comme le sous-groupe de matrices de la forme
$\pmatrix {1 &0 \cr 0& A}, A  \in SL(n, {\Bbb R})$ forment 
un sous-groupe isomorphe \`a $SL(n, {\Bbb Z})$ que nous notons
$\Gamma_0$. Les sous-groupes ab\'eliens unipotents 
constitu\'es des \'el\'ements de la forme
$\pmatrix {1 & * \cr 0& 1}$, et $\pmatrix {1 & 0 \cr * & 1}$
sont not\'es  $N^+$ et $N^-$ respectivement. Ils sont isomorphes
\`a ${\Bbb Z}^n$, et donnent lieu \`a deux produits semi-directs
$\Gamma_0 \ltimes N^+$ et $\Gamma_0 \ltimes N^-$ isomorphes
au produit semi-direct $SL(n, {\Bbb Z}) \ltimes {\Bbb Z}^n$ 
(o\`u 
$SL(n, {\Bbb Z})$ agit sur ${\Bbb Z}^n$ vu comme le r\'eseau 
canonique de 
${\Bbb R}^n$, dans un cas via la repr\'esentation canonique de 
$SL(n, {\Bbb R})$, et dans l'autre cas, via sa repr\'esentation duale).

Par hypoth\`ese, $\Gamma_0$ fixe un point, disons $x_0$, nous supposons que sa repr\'esentation
infinit\'esimale est la canonique. Maintenant, l'id\'ee est  de rev\^etir 
$M$ d'une structure projective $\Gamma$-invariante, en suivant la recette d\'etermin\'ee
par la structure projective usuelle (sur $S^n$ ou ${\Bbb R}P^n$).

Nous  montrons d'abord, que tout comme dans le cas standard, au moins \`a
indice fini pr\`es,  $N^+$ fixe $x_0$. L'id\'ee est que le centralisateur 
dans $\Gamma_0$, d'un \'el\'ement 
%quelconque 
$\gamma$ de $N^+$, est suffisamment grand, et ne peut avoir que des points fixes isol\'es, 
qui sont, dans leur ensemble, invariants par $\gamma$.

Ainsi, $\Gamma_0 \ltimes N^+$ admet $x_0$ comme  point fixe , et l'action 
de $\Gamma_0$  y est (localement) lin\'earisable.  Nous  montrons alors que 
$\Gamma_0 \ltimes N^+$ pr\'eserve (localement)  une structure projective plate  au
 voisinage 
de $x_0$ (celle d\'etermin\'ee par la connexion affine locale pr\'eserv\'ee 
par $\Gamma_0$).

Soit ${\cal S}(x_0, \Gamma_0)$ l'ouvert de Siegel de $x_0$ relatif \`a
l'action de $\Gamma_0$. Il correspond n\'ecessairement au cas dissipatif
du Th\'eor\`eme 1 (le cas conservatif est exclu, car il y a un groupe plus grand, $\Gamma$, 
qui agit). Il s'identifie donc \`a un ouvert de ${\Bbb R}^n$ invariant par 
l'action de $SL(n, {\Bbb Z})$ ($= \Gamma_0$). 

Le bord de ${\cal S}(x_0, \Gamma_0)$ dans $M$ peut \^etre 
tr\`es compliqu\'e.

La connexion projective s'\'etend, en particulier,  \`a ${\cal S}(x_0, \Gamma_0)$, mais aussi 
individuellement aux images $\rho(\gamma) ({\cal S}(x_0, \Gamma_0))$, $\gamma
\in N^+$ (car $N^+$ pr\'eserve la connexion projective).

Nous  montrons que la connexion projective plate s'\'etend  au voisinage de la fronti\`ere 
de ${\cal S}(x_0, \Gamma_0)$.  En d\'eveloppant la situation (\'equivariante) 
dans le substratum de la g\'eom\'etrie projective plate (i.e. 
$S^n$ ou ${\Bbb R}P^n$), nous identifierons le bord, et par suite l'adh\'erence 
 de ${\cal S}(x_0, \Gamma_0)$. Il
s'agit d'une h\'emisph\`ere (projective) dont le bord est soit $S^{n-1}$  soit 
${\Bbb R}P^{n-1}$. Dans ce dernier cas, l'adh\'erence de 
${\cal S}(x_0, \Gamma_0)$ couvre toute la vari\'et\'e 
$M$, qui sera donc ${\Bbb R}P^n$. Dans le cas o\`u  le bord est 
$S^{n-1}$, nous nous aidons de $N^-$ pour trouver une deuxi\`eme
h\'emisph\`ere compl\'ementaire \`a ${\cal S}(x_0, \Gamma_0)$, ce qui
permet d'identifier $M$ \`a
$S^n$.

\end{esqpreuve}

Il est probable que ce th\'eor\`eme 
soit vrai sans l'hypoth\`ese de compacit\'e de $M$
(l'\'enonc\'e serait alors qu'une telle action n'existe pas 
dans ce cas).

On en d\'eduit en particulier la rigidit\'e locale (parmi les 
actions analytiques), mais aussi une propri\'et\'e de 
non-d\'eg\'en\'erescence, comme dans le cas du tore ci-dessus.

\begin{corollaire} Soit $\Gamma$ un sous-groupe d'indice fini 
de $SL(n+1, {\Bbb Z})$, $n\geq 3$. Alors l'orbite par $Diff^\omega(S^n)$ de 
son action standard sur $S^n$ est ouverte-ferm\'ee dans 
$Rep(\Gamma,  Diff^\omega(S^n))$ muni de la topologie $C^1$ et elle
est  
ferm\'ee  au sens de la topologie $C^0$.
\end{corollaire}

%\paragraph{Cocycles.}

%\paragraph{Actions sur le tore.}

%{00}
\medskip
\noindent
CNRS, UMPA, \'Ecole Normale Sup\'erieure de Lyon \\
46, all\'ee d'Italie,
 69364 Lyon cedex 07,  FRANCE \\
Zeghib@umpa.ens-lyon.fr \\
http://umpa.ens-lyon.fr/\~{}zeghib/

\end{document}